\documentclass[12pt]{article}
\usepackage{amssymb,amsmath}
\usepackage{cases}
\usepackage{amsfonts}
\usepackage{color,xcolor}
\usepackage[left=2.0cm,right=2.0cm,top=2.0cm,bottom=2.0cm]{geometry}
\usepackage[colorlinks,citecolor=blue,urlcolor=blue]{hyperref}

\newtheorem{theorem}{Theorem}[section]

\newtheorem{remark}{Remark}[section]

\newcommand{\bal}{\begin{align}}
\newcommand{\bbal}{\begin{align*}}
\newcommand{\beq}{\begin{equation}}
\newcommand{\eeq}{\end{equation}}
\newcommand{\bca}{\begin{cases}}
\newcommand{\eca}{\end{cases}}
\def\div{\mathord{{\rm div}}}
\newcommand{\pa}{\partial}
\newcommand{\fr}{\frac}
\newcommand{\na}{\nabla}

\newcommand{\cd}{\cdot}

\newcommand{\dd}{\mathrm{d}}
\newcommand{\B}{\dot{B}}

\newcommand{\R}{\mathbb{R}}

\newcommand{\les}{\lesssim}

\begin{document}
\title{Global regularity for a class of 3D tropical climate model without thermal diffusion}

\author{Jinlu Li$^{1}$\footnote{E-mail: lijinlu@gnnu.cn}\quad  and Yanghai Yu$^{2}$\footnote{E-mail: yuyanghai214@sina.com( Corresponding author)}\\
\small $^1$\it School of Mathematics and Computer Sciences, Gannan Normal University, Ganzhou 341000, China\\
\small $^2$\it School of Mathematics and Statistics, Anhui Normal University, Wuhu, Anhui, 241002, China}

\date{\today}

\maketitle\noindent{\hrulefill}

{\bf Abstract:} In this paper, we establish the global regularity for the 3D tropical climate model with fractional dissipation.


\vskip0mm\noindent{\hrulefill}

\section{Introduction}\label{sec1}
This paper focuses on the following 3D tropical climate model given by
\begin{eqnarray}\label{3d}
        \left\{\begin{array}{ll}
          \partial_tu+u\cd\na u+\mu\Lambda^{2\alpha} u+\na p+\div(v\otimes v)=0,& x\in \R^3,t>0,\\
          \partial_tv+u\cd\na v+v\cd\na u+\na \theta+\nu\Lambda^{2\beta} v=0,& x\in \R^3,t>0,\\
          \partial_t\theta+u\cd\na \theta+ \eta\Lambda^{2\gamma}\theta+\div v=0,& x\in \R^3,t\geq0,\\
                   (u,v,\theta)|_{t=0}=(u_0,v_0,\theta_0),& x\in \R^3.\end{array}\right.
        \end{eqnarray}
here $\mu,\nu,\eta,\alpha,\beta,\gamma$ are non-negative parameters, $u=u(t,x)$ and $v=v(t,x)$ stand for the barotropic mode and the first baroclinic mode of the vector velocity, respectively. $\Pi=\Pi(t,x)$ and $\theta=\theta(t,x)$ denote the scalar pressure and scalar temperature, respectively. $\Lambda=\sqrt{-\Delta}$ is the Zygmund operator and the fractional power operator $\Lambda^{\delta}=(-\Delta)^{\frac{\delta}{2}}$ is defined through Fourier transform as
$$
\Lambda^{\delta} f(\xi)=\mathcal{F}^{-1}|\xi|^{\delta}\mathcal{F}f(\xi).
$$

The original version of (\ref{3d}) without any fractional Laplacian terms was derived by Frierson--Majda--Pauluis \cite{Frierson 2004} from the inviscid  primitive equations with the aid of performing a Galerkin truncation to the hydrostatic Boussinesq equations, of which the first baroclinic mode had been originally used in some studies of tropical atmosphere. More relevant background on the tropical climate model can be found in \cite{Gill 1980,Matsuno 1966,Biello 2003} and the references therein. From the mathematical point of view, the tropical climate model (\ref{3d}) are significant generalizations of the generalized magnetohydrodynamic (GMHD) equations which model the complex interaction between the fluid dynamic phenomena. In fact, when the temperature $\theta\equiv $ Constant, (\ref{3d}) is reduced to the GMHD equations, namely,
\begin{eqnarray}\label{Equ1.2}
        \left\{\begin{array}{ll}
          \partial_t u + u\cdot\nabla u +\nabla \Pi+\mu\Lambda^{2\alpha} u-v\cdot\nabla v=0,& x\in \mathbb{R}^d,t>0,\\
          \partial_t v + u\cdot\nabla v +\nu\Lambda^{2\beta} v-v\cdot\nabla u= 0,& x\in \mathbb{R}^d,t>0,\\
          {\rm{div}}u={\rm{div}}v=0,& x\in \mathbb{R}^d,t\geq0,\\
          (u(0,x),v(0,x))=(u_{0}(x),v_{0}(x)),& x\in \mathbb{R}^d,\end{array}\right.
        \end{eqnarray}
here $u=u(t,x)$ and $v=v(t,x)$ denote the velocity and the magnetic field of the fluid, respectively. The mathematical studies on the tropical climate model (\ref{3d}) have attracted considerable attention recently from various authors and have motivated a large number of research papers. Here, we mainly recall some global regularity results which are more relative with our research in this field. Li--Titi \cite{Li 2016} introduced a new quantity to bypass the obstacle caused by the absence of thermal diffusion and proved the global well-posedness of strong solutions for the 2D tropical climate model (\ref{3d}) with $\alpha=\beta=1$ and $\mu>0,\nu>0,\eta=0$. Later, Ye \cite{Ye 2016} obtained the global regularity of a tropical climate model with the very weak dissipation of the barotropic ($\alpha>0,\beta=\gamma=1$ and $\mu,\nu,\eta>0$) by the ``weakly nonlinear" energy estimates approach and the maximal $L_t^{q}L_x^{p}$ regularity for the heat kernel. Dong et al. \cite{Dong1 2018} established the global regularity results for the 2D system (\ref{3d}) without thermal diffusion with $\alpha+\beta=2,\beta\in(1,\fr32],\mu,\nu>0$ and $\alpha=2,\mu>0,\nu=\eta=0$, respectively. Recently, Dong et al. \cite{Dong2 2018} established the global existence and regularity for the 2D system (\ref{3d}) ($\mu,\nu,\eta>0,\beta=1$) with the fractional dissipation which are in two very broad ranges, namely,
$$\gamma\geq\fr{4+\alpha-\sqrt{\alpha^2+8\alpha+8}}{2},\quad\mbox{if}\quad0<\alpha<\fr12;\quad\gamma\geq1-\alpha,\quad\mbox{if}\quad\fr12\leq\alpha\leq1.$$

For the 3D case, Zhu \cite{Zhu 2018} obtained the global regularity for the tropical climate model (\ref{3d}) with only fractional diffusion on barotropic mode , namely, $\alpha\geq\fr52$, $\mu>0,\nu=\eta=0$. Naturally, we ask that how weak dissipation which the $u$-equations possess to ensure that the 3D system (\ref{3d}) has the global regularity? This is our principal goal of this paper. Our main result reads as follows

\begin{theorem}\label{the1.1} Consider the system \eqref{3d} with $\alpha=\beta=\fr54$ and $\gamma=0$. Assume that $(u_0, v_0,\theta_0)\in (H^2(\mathbb{R}^3)\cap \B^{0}_{2,1}(\mathbb{R}^3))^3$ and $\div u_0=0$. Then \eqref{3d}  has a unique global solution $(u, v,\theta)$ satisfying for any $T>0$
\begin{eqnarray*}
 (u, v,\theta)\in (L^{\infty}(0,T;H^{2}))^3\quad\mbox{and}\quad (u,v)\in L^{2}(0,T;H^{\fr{13}{4}}).
\end{eqnarray*}
\end{theorem}
\begin{remark}\label{rem1.0} As mentioned above, Theorem \ref{the1.1} holds true for the 3D GMHD system.
\end{remark}

\section{Proof of Theorem \ref{the1.1}}
Taking the $L^2$ inner product of Equs.$\eqref{3d}_1$--$\eqref{3d}_3$ with $u$, $v$ and $\theta$, respectively, then using the fact $\div u=0$, we have the following basic energy estimate
\bal\label{2.1}
||u||^2_{L^2}+||v||^2_{L^2}+||\theta||^2_{L^2}+2\int_0^T(||\Lambda^{\frac54}u||^2_{L^2}+||\Lambda^{\frac54}v||^2_{L^2})\dd s=||u_0||^2_{L^2}+||v_0||^2_{L^2}+||\theta_0||^2_{L^2}.
\end{align}
Applying $\dot{\Delta}_j$ to $\eqref{3d}_1$ and taking the $L^2$ inner product of the resulting equation with $\dot{\Delta}_j u$ yields
\bal\label{2.2}
\frac12\frac{\dd}{\dd t}||\dot{\Delta}_j u||_{L^2}+c2^{\fr52j}||\dot{\Delta}_j u||_{L^2}\leq ||\dot{\Delta}_j u_0||_{L^2}+||\dot{\Delta}_j\div(u\otimes u)||_{L^2}+||\dot{\Delta}_j\div(v\otimes v)||_{L^2},
\end{align}
which implies
\bal\label{2.3}
||u||_{L^\infty_T(\B^{0}_{2,1})}+||u||_{L^1_T(\B^{\frac52}_{2,1})}\lesssim ||u_0||_{\B^{0}_{2,1}}+\int^T_0(||u\otimes u||_{\B^{1}_{2,1}}+||v\otimes v||_{\B^{1}_{2,1}})\dd s.
\end{align}
Invoking the product estimate of two functions in homogeneous spaces (see Corollary 2.55 in \cite{Bahouri2011}):
\bal\label{2.4}
||f g||_{\B^{s+t}_{2,1}}\lesssim ||f||_{\dot{H}^{s}}||g||_{\dot{H}^{t}}\quad \mbox{with}\quad -3/2<s,t<3/2 \quad\mbox{and} \quad s+t>0
\end{align}
Therefore, it follows from \eqref{2.1} and \eqref{2.3}
\bal\label{2.5}
||u||_{L^\infty_T(\B^{0}_{2,1})}+||u||_{L^1_T(\B^{\frac52}_{2,1})}\lesssim ||u_0||_{\B^{0}_{2,1}}+\int^T_0(||u||^2_{\dot{H}^{\frac54}}+||v||^2_{\dot{H}^{\frac54}})\dd s\leq C.
\end{align}
Taking the $L^2$ inner product of Equs.$\eqref{3d}_1$--$\eqref{3d}_3$ with $-\Delta u$, $-\Delta v$ and $-\Delta \theta$, respectively, then using the fact $\div u=0$, we have
\bal\label{2.6}
&\frac12\frac{\dd}{\dd t}(||\na u||^2_{L^2}+||\na v||^2_{L^2}+||\na \theta||^2_{L^2})+||\Lambda^{\frac94}u||^2_{L^2}+||\Lambda^{\frac94}v||^2_{L^2}
\nonumber\\=&\underbrace{\int_{\R^3}(u\cd \na u)\cd \Delta u\dd x}_{I_1}+\underbrace{\int_{\R^3}\div(v\otimes v)\cd \Delta u\dd x}_{I_2}+\underbrace{\int_{\R^3}(v\cd\na u)\cd\Delta v\dd x}_{I_3}\nonumber\\
\quad&+\underbrace{\int_{\R^3}(u\cd\na v)\cd\Delta v\dd x}_{I_4}+\underbrace{\int_{\R^3}(u\cd\na \theta)\cd \Delta\theta \dd x}_{I_5},
\end{align}
Next, we need to estimate the above five terms one by one. \\
Notice that the fact $\div(v\otimes v)=v\cdot \na v+v\div v$, we have
\bbal
&I_1\lesssim ||u||_{L^6}||\na u||_{L^{\frac{12}{5}}}||\Delta u||_{L^{\frac{12}{5}}}\lesssim ||\na u||^2_{L^2}||\Lambda^{\frac54}u||^2_{L^2}+\varepsilon||\Lambda^{\frac94}u||^2_{L^2},
\\&I_2\lesssim ||v||_{L^6}||\na v||_{L^{\frac{12}{5}}}||\Delta u||_{L^{\frac{12}{5}}}\lesssim ||\na v||^2_{L^2}||\Lambda^{\frac54}v||^2_{L^2}+\varepsilon||\Lambda^{\frac94}u||^2_{L^2},
\\&I_3\lesssim ||v||_{L^6}||\na u||_{L^{\frac{12}{5}}}||\Delta v||_{L^{\frac{12}{5}}}\lesssim ||\na v||^2_{L^2}||\Lambda^{\frac54}u||^2_{L^2}+\varepsilon||\Lambda^{\frac94}v||^2_{L^2},
\\&I_4\lesssim ||u||_{L^6}||\na v||_{L^{\frac{12}{5}}}||\Delta v||_{L^{\frac{12}{5}}}\lesssim ||\na u||^2_{L^2}||\Lambda^{\frac54}v||^2_{L^2}+\varepsilon||\Lambda^{\frac94}v||^2_{L^2},\\
&I_5=\int_{\R^3}\pa_ku^i\pa_i\theta\pa_k\theta\dd x\les ||\na u||_{L^\infty}||\na \theta||^2_{L^2},
\end{align*}
where we have used the embedding $\dot{H}^\fr{1}{4}(\R^3)\hookrightarrow L^{\frac{12}{5}}(\R^3)$ and the Young inequality.

Putting the above estimates together with \eqref{2.6}, we obtain
\bbal
&\frac{\dd}{\dd t}||\na u,\na v,\na \theta||^2_{L^2}+||\Lambda^{\frac94}u,\Lambda^{\frac94}v||^2_{L^2}
\les||\na u,\na v,\na \theta||^2_{L^2}(||\Lambda^{\frac54}u,\Lambda^{\frac54}v||^2_{L^2}+||\na u||_{L^\infty}),
\end{align*}
combining the embedding $\dot{B}_{2,1}^{\fr32}(\R^3)\hookrightarrow L^\infty(\R^3)$, which implies
\bal\label{2.7}
\sup_{0\leq s\leq T}||\na u,\na v,\na \theta||^2_{L^2}+\int^T_0||\Lambda^{\frac94} u||^2_{L^2}+||\Lambda^{\frac94} v||^2_{L^2}\dd s\leq C.
\end{align}
Multiplying \eqref{2.2} by $2^j$ and summing over $j\in \mathbb{Z}$ yields
\bal\label{y}
||u||_{L^\infty_T(\B^{1}_{2,1})}+||u||_{L^1_T(\B^{\frac72}_{2,1})}\lesssim ||u_0||_{\B^{1}_{2,1}}+\int^T_0(||u\cdot \na u||_{\B^{1}_{2,1}}+||v\cdot \na v||_{\B^{1}_{2,1}}+||v\div v||_{\B^{1}_{2,1}})\dd s.
\end{align}
Using the product estimate \eqref{2.4} again, one has
\bbal
&||u\cdot \na u||_{\B^{1}_{2,1}}+||v\cdot \na v||_{\B^{1}_{2,1}}+||v\div v||_{\B^{1}_{2,1}}\lesssim ||u||_{\dot{H}^{\frac54}}||u||_{\dot{H}^{\frac94}}+||v||_{\dot{H}^{\frac54}}||v||_{\dot{H}^{\frac94}}.
\end{align*}
Therefore, we have from \eqref{y}
\bbal
||u||_{L^\infty_T(\B^{1}_{2,1})}+||u||_{L^1_T(\B^{\frac72}_{2,1})}\lesssim ||u_0||_{\B^{1}_{2,1}}+\int^T_0(||u||^2_{\dot{H}^{\frac54}}+||v||^2_{\dot{H}^{\frac54}}+||u||^2_{\dot{H}^{\frac94}}+||v||^2_{\dot{H}^{\frac94}})\dd s\leq C,
\end{align*}
which gives
\bal\label{2.9}
\int_0^T||\nabla^2 u||_{L^\infty}\dd s\leq C.
\end{align}
Applying $\Delta$ to both sides of $\eqref{3d}_1$, $\eqref{3d}_2$ and $\eqref{3d}_3$, taking the $L^2$ inner product of the resulting equations with $\Delta u$, $\Delta v $ and $\Delta \theta$, respectively, then adding them together yields that
\begin{align}\label{2.10}
&\frac{1}{2}\frac{\dd}{\dd t}||\Delta u,\Delta v,\Delta \theta||^2_{2}+||\Lambda^{\fr{13}{4}} u||_2^2+||\Lambda^{\fr{13}{4}}v||^2_{2}\nonumber\\
=&-\underbrace{\int_{\R^3}\Delta(u\cdot\nabla u)\cdot\Delta u\dd x}_{J_1}-\underbrace{\int_{\R^3}\Delta{\rm{div}}(v\otimes v)\cdot\Delta u\dd x}_{J_2}
-\underbrace{\int_{\R^3}\Delta(u\cdot\nabla v)\cdot\Delta v\dd x}_{J_3}\nonumber\\&-\underbrace{\int_{\R^3}\Delta(v\cdot\nabla u)\cdot\Delta v\dd x}_{J_4}-\underbrace{\int_{\R^3}\Delta(u\cdot\nabla\theta)\cdot\Delta \theta \dd x}_{J_5},
\end{align}
where we have used the fact
\begin{align*}
\int_{\R^3}\Delta\nabla\theta\cdot\Delta v\dd x+\int_{\R^3}\Delta{\rm{div}}v\cdot\Delta \theta \dd x=0.
\end{align*}
Next, we need to estimate the above four terms one by one as follows.

For the first term $J_1$, due to the divergence free condition $\div u=0$, one has
\begin{align}\label{2.11}
J_1=&~\int_{\R^3}\Delta u\cdot\nabla u\cdot\Delta u\dd x+\int_{\R^3}\partial_iu\cdot\partial_i\nabla u\cdot\Delta u\dd x
\les||\nabla u||_{L^\infty}||\Delta u||^{2}_{L^2}.
\end{align}
Similarly, we have
\begin{align}\label{2.12}
J_3=&\int_{\R^3}\Delta u\cdot\nabla v\cdot\Delta v\dd x+\int_{\R^3}\partial_iu\cdot\partial_i\nabla v\cdot\Delta v\dd x\nonumber\\
\les&~ ||\Delta u||_{L^\infty}||\nabla v||_{L^2}||\Delta v||_{L^2}+||\nabla u||_{L^\infty}||\Delta v||^{2}_{L^2},\\
J_5=&\int_{\R^3}\Delta u\cdot\nabla \theta\cdot\Delta \theta\dd x+\int_{\R^3}\partial_iu\cdot\partial_i\nabla \theta\cdot\Delta \theta\dd x\nonumber\\
\les&~||\Delta u||_{L^\infty}||\nabla \theta||_{L^2}||\Delta \theta||_{L^2}+ ||\nabla u||_{L^\infty}||\Delta \theta||^{2}_{L^2},
\end{align}
For the second term $J_2$, we decompose it into as
\bal\label{2.13}
J_2=&\underbrace{\int_{\R^3}\Delta (v\cdot \na v)\cdot\Delta u\dd x}_{J_{21}}+\underbrace{\int_{\R^3}\Delta (v\div v)\cdot\Delta u\dd x}_{J_{22}}.
\end{align}
We estimate $J_{21}$ as follows
\bbal
J_{21}\leq& \int_{\R^3}\Delta v\cdot \na v\cdot\Delta u\dd x+\int_{\R^3}\partial_i v\cdot \partial_i\na v\cdot\Delta u\dd x+\int_{\R^3} v\cdot \na \Delta v\cdot\Delta u\dd x\nonumber\\
\les&~||\Delta u||_{L^\infty}||\nabla v||_{L^2}||\Delta v||_{L^2}+||v||_{L^6}||\Delta \nabla v||_{L^\fr{12}{5}}||\Delta u||_{L^2}\nonumber\\
\les&~||\Delta u||_{L^\infty}(||\nabla v||^2_{L^2}+||\Delta v||^2_{L^2})+ ||\nabla v||^2_{L^2}||\Delta u||^2_{L^2}+\varepsilon||\Lambda^\fr{13}{4}v||^2_{L^2},
\end{align*}
Similarly, one has
\bbal
J_{22}\les&~||\Delta u||_{L^\infty}(||\nabla v||^2_{L^2}+||\Delta v||^2_{L^2})+ ||\nabla v||^2_{L^2}||\Delta u||^2_{L^2}+\varepsilon||\Lambda^\fr{13}{4}v||^2_{L^2},
\end{align*}
For the term $J_4$, we get
\bal\label{2.15}
J_{4}\leq& \int_{\R^3}\Delta v\cdot \na u\cdot\Delta v\dd x+\int_{\R^3}\partial_i v\cdot \partial_i\na u\cdot\Delta v\dd x+\int_{\R^3} v\cdot \na \Delta v\cdot\Delta v\dd x\nonumber\\
\les&~||\nabla u||_{L^\infty}||\Delta v||^2_{L^2}+||\nabla^2u||_{L^\infty}||\nabla v||_{L^2}||\Delta v||_{L^2}+||v||_{L^6}||\Delta \nabla v||_{L^\fr{12}{5}}||\Delta v||_{L^2}\nonumber\\
\les&~(||\nabla u||_{L^\infty}+||\nabla^2 u||_{L^\infty})(||\nabla v||^2_{L^2}+||\Delta v||^2_{L^2})+ ||\nabla v||^2_{L^2}||\Delta u||^2_{L^2}+\varepsilon||\Lambda^\fr{13}{4}v||^2_{L^2},
\end{align}
Combining the above estimates with \eqref{2.10}, one has
\begin{align*}
\frac{\dd}{\dd t}||\Delta u,\Delta v,\Delta \theta||^2_{L^2}+||\Lambda^{\fr{13}{4}} u,\Lambda^{\fr{13}{4}}v||_{L^2}^2
\leq C\Big(1+||\na u||_{\infty}+||\na^2 u||_{L^\infty}\Big)(1+||\Delta u,\Delta v,\Delta\theta||^2_{L^2}),
\end{align*}
which follows from \eqref{2.9} and Gronwall's Lemma that
\begin{align*}
\sup_{0\leq s\leq T}||\Delta u,\Delta v,\Delta \theta||^2_{L^2}+\int_0^T||\Lambda^{\fr{13}{4}}u,\Lambda^{\fr{13}{4}}v||_{L^2}^2\dd s\leq C(T).
\end{align*}
Thus, we complete the proof of Theorem \ref{the1.1}.
\section*{Acknowledgments} J. Li is supported by NSFC (No.11801090). Y. Yu is supported by NSF of Anhui Province (No.1908085QA05).

\end{document}